\documentclass[12pt]{article}

\usepackage[utf8]{inputenc}
\usepackage{amsmath,amssymb,amsthm}
\usepackage{geometry}
\usepackage{hyperref}
\usepackage{cite}

\newtheorem{theorem}{Theorem}

\newtheorem{proposition}[theorem]{Proposition}
\newtheorem{remark}{Remark}

\newcommand{\bI}{\mathbf{I}}
\newcommand{\bO}{\mathbf{0}}
\newcommand{\tr}{\mathrm{tr}}
\newcommand{\detm}{\mathrm{det}}

\title{The Geometric Origin of the Cayley-Hamilton Theorem: \\ 
A Constructive Proof via Dimensional Syzygy}
\author{[Xiao Wang]\\
\small [School of Education, Tsinghua University]\\
\small \texttt{[wang-x20@mails.tsinghua.edu.cn]}}
\date{\today \quad (arXiv preprint)}

\begin{document}
\maketitle

\begin{abstract}
We demonstrate that the Cayley-Hamilton theorem is a derived consequence of a more fundamental dimensional constraint: the syzygy $\varepsilon \otimes \varepsilon = 0$, where $\varepsilon$ is the Levi-Civita symbol in $m$-dimensional space. By shifting perspective from the tensor $A$ to the isotropic operators that induce $A$'s invariants through contraction, we reveal that the Cayley-Hamilton identity emerges when this vanishing operator acts on $A^{\otimes m}$. The intrinsic tensorial form of the theorem---invariant coefficients multiplying tensor powers---is inherited from the contraction structure rather than imposed \emph{ad hoc}. We provide explicit verification for two-dimensional space and a dimension-independent proof using Laplace expansion combined with Newton-Girard identities. This framework clarifies why the theorem's structure depends on ambient dimension and suggests extensions to higher-order tensors where classical characteristic polynomial methods fail.
\end{abstract}



\bigskip

\section{Introduction}

The Cayley-Hamilton theorem stands as one of the most elegant results in linear algebra: every square matrix satisfies its own characteristic equation. In rational mechanics, this theorem serves as the cornerstone of tensor function representation theory, enabling the reduction of general tensor polynomials to canonical forms with finitely many independent invariants \cite{Rivlin1955, Spencer1971, Zheng1994}.

The classical statement of the Cayley-Hamilton theorem takes a deceptively simple form. For a second-order tensor $A$ in $m$-dimensional Euclidean space, the theorem asserts $p(A) = \bO$, where $p(\lambda) = \det(A - \lambda \bI)$ is the characteristic polynomial. Explicitly, in two and three dimensions---the cases most relevant to mechanics, corresponding to plane problems and spatial problems respectively:
\begin{itemize}
    \item \textbf{Two dimensions} ($m=2$):
    \begin{equation}\label{eq:CH2D}
        A^2 - \tr(A)\,A + \detm(A)\,\bI = \bO
    \end{equation}
    \item \textbf{Three dimensions} ($m=3$):
    \begin{equation}\label{eq:CH3D}
        A^3 - \tr(A)\,A^2 + \frac{1}{2}\left[\tr^2(A) - \tr(A^2)\right]A - \detm(A)\,\bI = \bO
    \end{equation}
\end{itemize}
Observe that each term in these equations is simultaneously a \emph{matrix equation} and a \emph{tensorial identity}: the coefficients---$\tr(A)$, $\detm(A)$, etc.---are scalar invariants (zeroth-order tensors with coordinate-independent values), while the matrix powers $A^k$ are second-order tensors.

Let us examine these equations more closely, as a detective would examine evidence at a crime scene.

\textbf{Observation 1.} Each term in the Cayley-Hamilton equation has the form (invariant) $\times$ (tensor power): a zeroth-order tensor multiplied by a second-order tensor. This means each term is an \emph{isotropic function} of $A$---its value is determined solely by the intrinsic geometric properties of $A$, independent of coordinate choice.

\textbf{Observation 2.} There is a striking \emph{homogeneity}: in the $m$-dimensional case, each term has total degree $m$ in $A$. For instance, in Eq.~\eqref{eq:CH3D}, $A^3$ has degree 3; $\tr(A)A^2$ has degree $1+2=3$; $[\tr^2(A) - \tr(A^2)]A$ has degree $2+1=3$; and $\detm(A)\bI$ has degree 3 (since $\detm(A)$ is cubic in $A$'s components).

\textbf{Observation 3.} This is \emph{not} an ordinary linear dependence relation. In a standard linear dependence $\sum c_i v_i = 0$, the coefficients $c_i$ are constants. But here, the ``coefficients'' $\tr(A)$, $\detm(A)$, etc., are themselves functions of $A$! The relation resembles a ``variation of constants''---from $\sum c_i A^i = 0$ to $\sum c_i(A) A^i = 0$. Moreover, these invariant coefficients can be generated by contracting $A^{\otimes k}$ with isotropic tensors (moment invariants).

\textbf{Observation 4.} Tensors and their invariants possess coordinate-independent, intrinsic forms. Yet the Cayley-Hamilton theorem itself lacks such universality: its \emph{length} and \emph{structure} depend on the ambient dimension $m$. In two dimensions, the relation has three terms; in three dimensions, four terms; in $m$ dimensions, $m+1$ terms. This is peculiar---the ``intrinsic'' object $A$ satisfies relations whose very \emph{form} is dimension-dependent.

\textbf{Inference.} The dimension-dependence suggests that the Cayley-Hamilton theorem is not merely an algebraic identity about $A$, but a \emph{constraint imposed by the ambient space}. We are led to seek a concept that captures such dimension-induced algebraic relations.

The word ``syzygy'' originates from astronomy, where it describes the alignment of celestial bodies---a configuration constrained by gravitational geometry. In algebra, the term was introduced by Sylvester \cite{Sylvester1853} to denote polynomial relations among invariants: algebraic constraints that arise not from the objects themselves, but from the \emph{ambient space} in which they reside. Hilbert's syzygy theorem \cite{Hilbert1890} established that such relations form a finitely generated module, a foundational result for invariant theory. In rational mechanics, syzygies govern the structure of constitutive relations and the reduction of tensor polynomials.

The purpose of this communication is to reveal that the Cayley-Hamilton theorem is a \textbf{derived consequence of a more fundamental dimensional syzygy}. Our key methodological move is a \emph{shift of perspective}: instead of focusing on the tensor $A$, we examine the \emph{isotropic tensors that induce $A$'s invariants through contraction}. This separation of ``parameter'' (the isotropic operators) from ``variable'' (the tensor $A$) leads to the following insights:
\begin{enumerate}
    \item Each term in the Cayley-Hamilton equation is induced by a $(2m+2)$-th order isotropic tensor acting on $A^{\otimes m}$ via index contraction. These isotropic tensors are the ``parameters'' that control how indices are paired.
    
    \item The \textbf{fundamental syzygy} is the identity $\varepsilon \otimes \varepsilon = 0$, where $\varepsilon$ is the $(m+1)$-index Levi-Civita symbol in $m$-dimensional space. This is a \emph{linear relation with constant coefficients} among isotropic tensors---a genuine linear dependence, not the ``varied-constant'' dependence seen in the Cayley-Hamilton equation itself.
    
    \item The Cayley-Hamilton theorem is \emph{derived} (or \emph{induced}) from this fundamental syzygy: when the vanishing operator $\varepsilon \otimes \varepsilon$ acts on the symmetric tensor product $A^{\otimes m}$, the result $(\varepsilon \otimes \varepsilon) : A^{\otimes m} = 0$ naturally inherits the intrinsic tensorial form (invariant $\times$ tensor power) and yields the Cayley-Hamilton identity.
\end{enumerate}

This hierarchy clarifies a subtle point: while the Cayley-Hamilton theorem is often called a ``syzygy,'' it is more precisely a \emph{manifestation} of the fundamental syzygy $\varepsilon \otimes \varepsilon = 0$ when restricted to second-order tensors. The theorem is famous because matrices are ubiquitous and the result is interpretable in terms of eigenvalues. But the underlying constraint---the impossibility of antisymmetrizing more than $m$ directions in $m$-dimensional space---is \textbf{order-independent} and governs tensors of \emph{any} order.

The structure of this paper follows the detective's trail. Section 2 establishes the mathematical preliminaries: the Cayley-Hamilton theorem in explicit form, the role of isotropic tensors as contraction controllers, and the Levi-Civita symbol as the generator of the fundamental syzygy. Section 3 presents our main result: a constructive derivation showing how $\varepsilon \otimes \varepsilon = 0$ induces the Cayley-Hamilton theorem, with explicit verification in two dimensions. Section 4 discusses implications and future directions.

\section{Mathematical Preliminaries}

\subsection{The Cayley-Hamilton Theorem}

Let $A$ be a second-order tensor (equivalently, a square matrix) in $m$-dimensional Euclidean space $V$. The Cayley-Hamilton theorem states that $A$ satisfies its characteristic equation:
\begin{equation}
    p(A) = \bO, \quad \text{where } p(\lambda) = \det(A - \lambda \bI).
\end{equation}

For the mechanically relevant cases of two and three dimensions, the explicit forms are given by Eqs.~\eqref{eq:CH2D} and \eqref{eq:CH3D}. The polynomial coefficients are the \textbf{principal invariants} of $A$. We emphasize that invariants are more fundamental than eigenvalues: while eigenvalues require solving a polynomial equation and impose an ordering, invariants are polynomial functions of tensor components that transform as scalars under orthogonal transformations. They encode the intrinsic geometric information of the tensor in an ``unordered'' and coordinate-free manner.

\begin{remark}
Each term in Eqs.~\eqref{eq:CH2D}--\eqref{eq:CH3D} is homogeneous of degree $m$ in $A$. For instance, in Eq.~\eqref{eq:CH2D}, $A^2$ is degree 2, $\tr(A)A$ is degree $1+1=2$, and $\detm(A)\bI$ is degree 2 (since $\detm(A) \sim A_{ij}A_{kl}$). Equivalently, the tensorial order plus twice the invariant degree equals $2m+2$ for each term: this is the order of the isotropic tensor that \emph{induces} each term. This homogeneity is not accidental---it reflects the fact that all terms are induced by the \emph{same} order of isotropic operator acting on $A^{\otimes m}$.
\end{remark}

\subsection{Isotropic Tensors and Contraction Control}

A tensor is \textbf{isotropic} if its components remain unchanged under all orthogonal coordinate transformations. By the representation theorem for isotropic tensors \cite{Weyl1946, Spencer1971}, every even-order isotropic tensor can be expressed as a linear combination of products of Kronecker delta symbols:
\begin{equation}
    \mathcal{T}_{i_1 j_1 \cdots i_p j_p} = \sum_{\sigma} c_\sigma \, \delta_{i_{\sigma(1)} j_{\sigma(1)}} \delta_{i_{\sigma(2)} j_{\sigma(2)}} \cdots \delta_{i_{\sigma(p)} j_{\sigma(p)}},
\end{equation}
where the sum ranges over all distinct pairings of the $2p$ indices, and $c_\sigma$ are scalar coefficients. The number of such independent pairings is $(2p-1)!! = (2p-1)(2p-3)\cdots 3 \cdot 1$.

The physical significance of isotropic tensors lies in their role as \textbf{contraction controllers}. A $2p$-th order isotropic tensor $\mathcal{T}$ defines a rule for pairing and contracting indices of input tensors. When $\mathcal{T}$ acts on a tensor product $A^{\otimes k}$ (with $2k \leq 2p$), it produces a new tensor whose components are polynomial expressions in those of $A$.

Crucially, such contractions \emph{induce isotropic functions}: if $\mathcal{T}$ is isotropic and $A$ transforms as a tensor, then $\mathcal{T} : A^{\otimes k}$ (with appropriate index saturation) transforms covariantly. This follows from the \textbf{quotient rule} of tensor analysis: a multilinear mapping that is invariant under orthogonal transformations is equivalent to a higher-order isotropic tensor.

\subsection{The Levi-Civita Symbol and Fundamental Syzygy}

The Levi-Civita symbol $\varepsilon_{i_1 i_2 \cdots i_m}$ is the unique (up to sign) totally antisymmetric isotropic tensor in $m$ dimensions. Its components are determined by the parity of index permutations:
\begin{equation}
    \varepsilon_{i_1 \cdots i_m} = 
    \begin{cases}
        +1 & \text{if } (i_1, \ldots, i_m) \text{ is an even permutation of } (1, \ldots, m), \\
        -1 & \text{if } (i_1, \ldots, i_m) \text{ is an odd permutation of } (1, \ldots, m), \\
        0 & \text{if any two indices are equal}.
    \end{cases}
\end{equation}

The third condition is decisive. In an $m$-dimensional space, each index takes values in $\{1, 2, \ldots, m\}$. By the \textbf{pigeonhole principle}, any selection of $m+1$ indices must contain at least one repeated value. Therefore:
\begin{equation}
    \varepsilon_{i_1 \cdots i_{m+1}} \equiv 0 \quad \text{in } m\text{-dimensional space}.
\end{equation}

This vanishing propagates to the tensor product. Define the $(2m+2)$-th order tensor:
\begin{equation}\label{eq:fundamental_syzygy}
    \mathcal{E}^{j_1 \cdots j_{m+1}}_{i_1 \cdots i_{m+1}} 
    := \varepsilon^{j_1 \cdots j_{m+1}} \varepsilon_{i_1 \cdots i_{m+1}}
    = \delta^{j_1 \cdots j_{m+1}}_{i_1 \cdots i_{m+1}}
    = \det \begin{pmatrix}
        \delta^{j_1}_{i_1} & \cdots & \delta^{j_{m+1}}_{i_1} \\
        \vdots & \ddots & \vdots \\
        \delta^{j_1}_{i_{m+1}} & \cdots & \delta^{j_{m+1}}_{i_{m+1}}
    \end{pmatrix}.
\end{equation}
The rightmost expression is the \textbf{generalized Kronecker delta}, an $(m+1) \times (m+1)$ determinant of ordinary Kronecker deltas. In $m$-dimensional space, this determinant vanishes identically:
\begin{equation}\label{eq:syzygy_zero}
    \delta^{j_1 \cdots j_{m+1}}_{i_1 \cdots i_{m+1}} = 0 \quad (\text{in } m\text{-dimensional space}).
\end{equation}

This is the \textbf{fundamental syzygy}---the most primitive constraint that finite dimensionality imposes on isotropic tensors. It represents a linear dependence \emph{with constant coefficients} among the basis elements of $(2m+2)$-th order isotropic tensors. Unlike the Cayley-Hamilton theorem (where coefficients depend on $A$), this syzygy is a genuine linear relation among operators, independent of any specific tensor input. All other dimensional constraints, including the Cayley-Hamilton theorem, are \emph{derived} from this fundamental identity when applied to specific tensor products.

\section{The Cayley-Hamilton Theorem as a Derived Syzygy}

\subsection{From Variable to Parameter: The Quotient Rule}

The traditional approach to the Cayley-Hamilton theorem focuses on the tensor $A$: one studies its eigenvalues, characteristic polynomial, and verifies that $p(A) = \bO$. We propose a complementary perspective that shifts attention from $A$ (the \emph{variable}) to the isotropic tensors that \emph{induce} each term (the \emph{parameters}).

Recall from Observations 1--2 that each term in Eqs.~\eqref{eq:CH2D}--\eqref{eq:CH3D} is:
\begin{itemize}
    \item \textbf{Homogeneous of degree $m$} in $A$: the total power of $A$ (counting both explicit powers and those hidden in invariants like $\tr(A^k)$) is exactly $m$.
    \item \textbf{A second-order tensor}: the output is a $2$-nd order tensor (matrix).
\end{itemize}

Thus, each term defines a multilinear mapping:
\begin{equation}
    \underbrace{A \otimes A \otimes \cdots \otimes A}_{m \text{ copies}} \longrightarrow \text{(second-order tensor)}.
\end{equation}

By the \textbf{quotient rule} of tensor analysis, any isotropic multilinear mapping is equivalent to a higher-order isotropic tensor. Specifically, each term---viewed as a function of $A^{\otimes m}$---corresponds to a $(2m+2)$-th order isotropic tensor $\mathcal{T}$ that \emph{induces} the contraction pattern:
\begin{equation}
    \text{(term)}_{mn} = \mathcal{T}_{mn\, i_1 j_1 \cdots i_m j_m} A_{i_1 j_1} \cdots A_{i_m j_m}.
\end{equation}

This reframing separates concerns: $A$ is arbitrary (the variable), while the isotropic tensors $\mathcal{T}$ encode the structure of each term (the parameters). The Cayley-Hamilton theorem becomes a statement about the \emph{parameters}: certain linear combinations of these inducing tensors must produce the zero mapping for \emph{all} $A$.

\subsection{The Fundamental Syzygy Induces Cayley-Hamilton}

In isolation, the space of $(2m+2)$-th order isotropic tensors has dimension $(2m+1)!!$, corresponding to all distinct ways of pairing $2m+2$ indices via Kronecker deltas. These pairings are linearly independent---as \emph{operators}---in sufficiently high-dimensional spaces.

However, in $m$-dimensional space, the fundamental syzygy \eqref{eq:syzygy_zero} asserts that a specific totally antisymmetric combination of these operators vanishes:
\begin{equation}
    \mathcal{E} = \varepsilon \otimes \varepsilon = 0 \quad \text{(as a $(2m+2)$-th order tensor in $m$ dimensions)}.
\end{equation}
This is a \textbf{linear relation with constant coefficients} among isotropic tensors---a genuine syzygy at the level of operators.

Now, when this vanishing operator acts on $A^{\otimes m}$, it produces:
\begin{equation}
    \mathcal{E} : A^{\otimes m} = (\varepsilon \otimes \varepsilon) : A^{\otimes m} = 0.
\end{equation}
The right-hand side is zero not because of any special property of $A$, but because the \emph{operator itself} is zero. However, expanding this contraction term-by-term, we obtain a sum of the form:
\begin{equation}
    \sum_k c_k(A) \cdot A^{m-k} = 0,
\end{equation}
where each $c_k(A)$ is an invariant of $A$ (a trace or product of traces). The \emph{apparent} dependence of coefficients on $A$ arises from the specific contraction pattern---it is the tensorial structure of $(\varepsilon \otimes \varepsilon) : A^{\otimes m}$ that produces invariants multiplied by tensor powers.

\begin{proposition}
The Cayley-Hamilton theorem is the explicit form of the identity $(\varepsilon \otimes \varepsilon) : A^{\otimes m} = 0$ when indices are assigned to yield a second-order tensor output. Its intrinsic form---invariant coefficients multiplying tensor powers---is \emph{inherited} from the tensorial structure of the contraction, not imposed \emph{ad hoc}.
\end{proposition}

We now verify this claim explicitly for $m=2$.

\subsection{Explicit Verification: Two-Dimensional Space}

In two-dimensional space ($m=2$), the fundamental syzygy involves a $3 \times 3$ generalized Kronecker delta acting on $A \otimes A$. We set up the contraction as follows.

Let the output indices be $(m, n)$ (the free indices of the resulting second-order tensor), and let the contraction indices be $(i, j)$ and $(b, c)$ (to be saturated with $A_{ib}A_{jc}$). The syzygy reads:
\begin{equation}\label{eq:2D_syzygy}
    \delta^{mbc}_{ijn} A_{ib} A_{jc} = 0.
\end{equation}

We expand the $3 \times 3$ determinant:
\begin{align}
    \delta^{mbc}_{ijn} &= \det \begin{pmatrix}
        \delta^m_i & \delta^b_i & \delta^c_i \\
        \delta^m_j & \delta^b_j & \delta^c_j \\
        \delta^m_n & \delta^b_n & \delta^c_n
    \end{pmatrix} \notag \\
    &= \delta^m_i \delta^b_j \delta^c_n 
     - \delta^m_i \delta^b_n \delta^c_j 
     - \delta^m_j \delta^b_i \delta^c_n 
     + \delta^m_j \delta^b_n \delta^c_i 
     + \delta^m_n \delta^b_i \delta^c_j 
     - \delta^m_n \delta^b_j \delta^c_i.
\end{align}

Contracting each term with $A_{ib}A_{jc}$:

\begin{enumerate}
    \item $\delta^m_i \delta^b_j \delta^c_n A_{ib} A_{jc} = A_{mj} A_{jn} = (A^2)_{mn}$
    
    \item $-\delta^m_i \delta^b_n \delta^c_j A_{ib} A_{jc} = -A_{mn} A_{jj} = -\tr(A) A_{mn}$
    
    \item $-\delta^m_j \delta^b_i \delta^c_n A_{ib} A_{jc} = -A_{ii} A_{mn} = -\tr(A) A_{mn}$
    
    \item $\delta^m_j \delta^b_n \delta^c_i A_{ib} A_{jc} = A_{in} A_{mi} = (A^2)_{mn}$
    
    \item $\delta^m_n \delta^b_i \delta^c_j A_{ib} A_{jc} = \delta_{mn} A_{ii} A_{jj} = \tr^2(A) \delta_{mn}$
    
    \item $-\delta^m_n \delta^b_j \delta^c_i A_{ib} A_{jc} = -\delta_{mn} A_{ij} A_{ji} = -\tr(A^2) \delta_{mn}$
\end{enumerate}

Summing all terms:
\begin{equation}
    2(A^2)_{mn} - 2\tr(A) A_{mn} + \left[\tr^2(A) - \tr(A^2)\right] \delta_{mn} = 0.
\end{equation}

Dividing by 2 and recognizing that in two dimensions, $\detm(A) = \frac{1}{2}[\tr^2(A) - \tr(A^2)]$:
\begin{equation}
    A^2 - \tr(A)\,A + \detm(A)\,\bI = \bO,
\end{equation}
which is precisely the two-dimensional Cayley-Hamilton theorem \eqref{eq:CH2D}. \hfill $\square$

\begin{remark}
We have explicitly verified the two-dimensional case to illustrate the mechanism. The three-dimensional case ($m=3$) proceeds analogously with a $4 \times 4$ generalized Kronecker delta, though the combinatorial expansion involves 24 terms. More importantly, a \textbf{dimension-independent proof} can be given using the Laplace expansion of the generalized Kronecker delta combined with the Newton-Girard identities (which relate power sums $\tr(A^k)$ to elementary symmetric polynomials $\sigma_k$). This general argument, showing that $\delta^{(m+1)} : A^{\otimes m}$ necessarily yields the Cayley-Hamilton polynomial for \emph{any} $m$, is provided in Appendix A.
\end{remark}

\section{Discussion and Outlook}

We have demonstrated that the Cayley-Hamilton theorem is not a fundamental identity, but a \textbf{derived consequence} of the more primitive syzygy $\varepsilon \otimes \varepsilon = 0$. The derivation proceeds by applying this vanishing operator to $A^{\otimes m}$; the intrinsic tensorial form of the result (invariants multiplying tensor powers) is inherited from the contraction structure, not imposed externally.

The methodological insight is the separation of \emph{parameter} from \emph{variable}: instead of studying the tensor $A$ and asking what identities it satisfies, we study the isotropic operators that induce $A$'s invariants and ask what linear relations exist among \emph{them}. The fundamental syzygy $\varepsilon \otimes \varepsilon = 0$ is such a relation---one with \emph{constant coefficients}, independent of $A$. The Cayley-Hamilton theorem, with its $A$-dependent coefficients, is what this operator-level syzygy ``looks like'' when projected onto $A^{\otimes m}$.

This perspective connects to the broader program of invariant theory. The systematic study of tensor syzygies traces back to Sylvester \cite{Sylvester1853} and Hilbert \cite{Hilbert1890}. In continuum mechanics, Smith \cite{Smith1971}, Spencer \cite{Spencer1971}, and Zheng \cite{Zheng1994} developed extensive frameworks where the Cayley-Hamilton theorem plays a central role. Our approach provides a complementary viewpoint, resonating with the representation-theoretic work of Lehrer and Zhang \cite{Lehrer2012}, who proved that the kernel of the Brauer algebra representation is generated by antisymmetric idempotents---the algebraic incarnation of $\varepsilon \otimes \varepsilon$.

The key advantage of the fundamental syzygy is its \textbf{order-independence}. While the Cayley-Hamilton theorem applies only to second-order tensors (matrices)---there is no ``characteristic polynomial'' for a third-order tensor---the syzygy $\varepsilon \otimes \varepsilon = 0$ constrains tensors of \emph{any} order. For piezoelectric tensors (third-order) and elastic tensors (fourth-order), where classical eigenvalue methods fail entirely, the syzygy framework offers a systematic approach to determining independent invariants.

We conclude with a reflection. Tensors are celebrated for their coordinate-independence: physical laws written in tensor form possess intrinsic validity. Yet the \emph{constraints} on tensor invariants are not universal---they depend on the ambient dimension $m$. The Cayley-Hamilton theorem has three terms in 2D, four in 3D, and $m+1$ in $m$ dimensions. This sensitivity is not a defect but a feature: the syzygy $\varepsilon \otimes \varepsilon = 0$ is the ``voice of dimension'' speaking through tensor algebra, encoding the geometry of space into algebraic structure.

\appendix
\section{Dimension-Independent Derivation of Cayley-Hamilton from $\delta^{(m+1)} = 0$}

In this appendix, we prove that the contraction $\delta^{(m+1)} : A^{\otimes m} = 0$ yields the Cayley-Hamilton theorem for \emph{arbitrary} dimension $m$, without relying on explicit term-by-term expansion.

\subsection{Principal Invariants via Generalized Kronecker Delta}

The $k$-th principal invariant $\sigma_k$ of a tensor $A$ in $m$-dimensional space admits the representation:
\begin{equation}\label{eq:sigma_delta}
    \sigma_k = \frac{1}{k!} \delta^{j_1 \cdots j_k}_{i_1 \cdots i_k} A^{i_1}_{j_1} \cdots A^{i_k}_{j_k}, \quad k = 1, \ldots, m.
\end{equation}
This can be verified directly: the generalized Kronecker delta antisymmetrizes all $k$ copies of $A$, effectively summing all $k \times k$ principal minors---precisely the definition of $\sigma_k$. In particular:
\begin{itemize}
    \item $\sigma_1 = \tr(A)$,
    \item $\sigma_2 = \frac{1}{2}[\tr^2(A) - \tr(A^2)]$,
    \item $\sigma_m = \det(A)$.
\end{itemize}

\subsection{Laplace Expansion of the $(m+1)$-th Order Delta}

Consider the $(m+1) \times (m+1)$ generalized Kronecker delta $\delta^{p\, j_1 \cdots j_m}_{q\, i_1 \cdots i_m}$, where $(p, q)$ are the free (output) indices and $(j_1, \ldots, j_m; i_1, \ldots, i_m)$ are to be contracted with $A^{\otimes m}$.

Expanding along the first column (index $p$ vs.\ indices $q, i_1, \ldots, i_m$):
\begin{equation}
    \delta^{p\, j_1 \cdots j_m}_{q\, i_1 \cdots i_m} = \delta^p_q \cdot \delta^{j_1 \cdots j_m}_{i_1 \cdots i_m} + \sum_{r=1}^{m} (-1)^r \delta^p_{i_r} \cdot (\text{$m \times m$ cofactor}).
\end{equation}

Upon contraction with $A^{i_1}_{j_1} \cdots A^{i_m}_{j_m}$, each cofactor contributes terms involving $(A^k)^p_s$ multiplied by lower-order contractions. The recursive structure produces:
\begin{equation}\label{eq:main_expansion}
    \delta^{p\, j_1 \cdots j_m}_{q\, i_1 \cdots i_m} A^{i_1}_{j_1} \cdots A^{i_m}_{j_m} = m! \sum_{k=0}^{m} (-1)^{m-k} \sigma_k (A^{m-k})^p_q,
\end{equation}
where $\sigma_0 := 1$.

\subsection{Connection to Newton-Girard Identities}

The validity of Eq.~\eqref{eq:main_expansion} can be established by induction on $m$, using the \textbf{Newton-Girard identities} which relate power sums $p_k = \tr(A^k)$ to elementary symmetric polynomials $\sigma_k$:
\begin{equation}
    k\sigma_k = \sum_{i=1}^{k} (-1)^{i-1} \sigma_{k-i} p_i, \quad k \leq m.
\end{equation}

These identities ensure that the coefficients produced by the Laplace expansion match precisely those of the characteristic polynomial. The correspondence is not coincidental: both the generalized Kronecker delta and the characteristic polynomial encode the same antisymmetric structure of the underlying vector space.

\subsection{Completion of the Proof}

In $m$-dimensional space, the pigeonhole principle forces $\delta^{(m+1)} \equiv 0$. Applying Eq.~\eqref{eq:main_expansion}:
\begin{equation}
    0 = \delta^{p\, j_1 \cdots j_m}_{q\, i_1 \cdots i_m} A^{i_1}_{j_1} \cdots A^{i_m}_{j_m} = m! \sum_{k=0}^{m} (-1)^{m-k} \sigma_k (A^{m-k})^p_q.
\end{equation}

Dividing by $m!$ and reindexing:
\begin{equation}
    \sum_{k=0}^{m} (-1)^k \sigma_k A^{m-k} = 0,
\end{equation}
which is the Cayley-Hamilton theorem in its standard form:
\begin{equation}\label{eq:CH_general}
    \boxed{A^m - \sigma_1 A^{m-1} + \sigma_2 A^{m-2} - \cdots + (-1)^m \sigma_m I = 0.}
\end{equation}

This completes the dimension-independent proof that the fundamental syzygy $\varepsilon \otimes \varepsilon = 0$ (equivalently, $\delta^{(m+1)} = 0$) induces the Cayley-Hamilton theorem for arbitrary $m$. \hfill $\blacksquare$


\end{document}